\documentclass[a4,11pt]{article}
\usepackage{amssymb,amsfonts,amsthm,amsmath}
\usepackage[brazil]{babel}
\usepackage{epsfig}
\usepackage[utf8]{inputenc}
\usepackage{mathrsfs}
\usepackage{graphicx}
\usepackage{amssymb}
\usepackage{pgfplots}
\usepackage{float}

\newtheorem{definition}{Definition}
\newtheorem{application}{Application}

\newtheorem{property}{Property}

\begin{document}

\title{ \Large{\bf Solutions of time-fractional differential equations using homotopy analysis method}}


\author{ {\bf {\large{D. S. Oliveira }}} \\
 {\small Coordination of Civil Engineering} \\
{\small UTFPR} \\
 {\small 85053-525, Guarapuava, PR, Brazil} \\
 {\small oliveiradaniela@utfpr.edu.br} \\ \\
 \and
  {\bf {\large{E. Capelas de Oliveira}}} \\
 {\small Department of Applied Mathematics} \\
{\small Imecc - Unicamp} \\
 {\small 13083-859, Campinas, SP, Brazil} \\
 {\small capelas@ime.unicamp.br} }

 \date{}

\maketitle

\thispagestyle{empty}
\vspace{-.4cm}


\noindent{\bf Abstract:} {We have used the homotopy analysis method to obtain solutions of linear and nonlinear 
fractional partial differential differential equations with initial conditions. We 
replace the first order time derivative by $\psi$-Caputo fractional derivative, and also we compare the 
results obtained by the homotopy analysis method with the exact solutions.}\\


\vspace*{.2cm}
\section{Introduction} 

Fractional order differential equations, systems of fractional algebraic-dif\-fer\-en\-tial equations and
fractional integrodifferential equations have been widely studied. Methods for obtaining analytical solutions 
to these problems, in its nonlinear form, are commonly used, among them the Adomian decomposition method 
\cite{Adomian,jafari2006,momani}, homotopy perturbation method (HPM) \cite{momani} and homotopy analysis method (HAM) 
\cite{Liao2,Liao}. Jafari and Seifi \cite{jafari2009} have obtained solutions for linear and nonlinear
fractional diffusion and wave equations by means of the HAM. Ganjiani \cite{ganjiani} has
discussed only nonlinear fractional differential equations and Xu et al. \cite{Xu} have discussed
fractional partial differential equations subject to the boundary conditions and initial condition,
both by means of homotopy analysis method. S{\l}ota et al. \cite{slota}
have applied HAM for solving integrodifferential equations and Zhang et al. \cite{zhang}
have investigate numerical solutions of higher-order fractional integrodifferential equations with
boundary conditions. Zurigat et al. \cite{zurigat2010,zurigat} have used HAM to solve systems of 
fractional algebraic-differential equations.

Since, we are interested in the analytical solution of fractional partial differential equations, we must choose a 
particular fractional derivative. There are several types of fractional derivatives defined in terms of a respective fractional 
integral \cite{almeida,danieco,graziane,vanterler}. Perhaps the various ways of approaching the fractional derivative 
reside in the fact that, until now, we don't have a classic geometric interpretation, as in the case of an 
ordinary differential where we associate the concept of derivative with the concept of tangency. 
Here, we choose the $\psi$-Caputo fractional derivative \cite{almeida} to discuss our applications. 
The choice of this fractional differentiation operator is due to the fact that when we derive a constant, the result is identically zero and as a particular case recovers the classical Caputo fractional derivative.

The paper is organized in the following way. In Section \ref{sec:2}, we present some basic definitions aof the 
fractional calculus. In Section \ref{sec:3}, we describe the HAM, and three examples are present in \ref{sec:4}.
The first approach, we discuss linear time-fractional diffusion equation;  
the second approach a nonlinear time-fractional gas-dynamic equation is investigated. The third approach, we discuss 
nonlinear time-fractional KdV equation. At the end of each application numerical solutions to show the 
efficiency of the method are presented. Concluding remarks close the paper.

\section{Fractional calculus}
\label{sec:2}

In this section we present some concepts of the fractional calculus that are useful in the remainder of the text.
\begin{definition}
Let $\alpha>0$, $I=[a,b]$ be a finite or infinite interval, $f$ an integrable function defined on $I$ and 
$\psi\in{C}^1(I)$ an increasing function such that $\psi'(x)\neq{0}$, for all $x\in I$. The left fractional
integral of $f$ with respect to another function $\psi$ of order $\alpha$ is defined as
\textnormal{\cite{almeida,kilbas}}
\begin{eqnarray}
I^{\alpha,\psi}_{a+}[f(x,t)]=\frac{1}{\Gamma(\alpha)}\int_{a}^{t}\psi'(\tau)(\psi(t)-\psi(\tau))^{\alpha-1}f(x,\tau)\textnormal{d}\tau. \label{integral}
\end{eqnarray}
For $\alpha=0$, we have
$$I^{0,\psi}_{a+}[f(x,t)]=f(x,t).$$
\end{definition}
\begin{definition}
Let $\alpha>0$, $n\in\mathbb{N}$, $I$ is the interval $-\infty\leq{a}<b\leq{\infty}$, $f,\psi\in{C^n}(I)$ two 
functions such that $\psi$ is increasing and $\psi'(x)\neq{0}$, for all $x\in I$. The left $\psi$-Caputo 
fractional derivative of $f$ of order $\alpha$ is given by \textnormal{\cite{almeida}}
$${^{C}D_{a+}^{\alpha,\psi}}[f(x,t)]=I_{a+}^{n-\alpha,\psi}
\left(\frac{1}{\psi'(t)}\frac{\partial}{\partial t}\right)^n f(x,t),$$
where
$$n=[\alpha]+1 \quad \mbox{for}\quad \alpha\notin\mathbb{N}, \quad\quad 
n=\alpha\quad \mbox{for} \quad \alpha\in\mathbb{N}.$$
To simplify notation, we will use the abbreviated notation
$$f^{[n],\psi}(x,t)=\left(\frac{1}{\psi'(t)}\frac{\partial}{\partial t}\right)^n f(x,t).$$
\end{definition}
\begin{property}
Let $f\in C^{n}[a,b]$, $\alpha>0$ and $\delta>0,$ \textnormal{\cite{almeida}}.
\begin{enumerate}

\item $f(t)=(\psi(t)-\psi(a))^{\delta-1}$, then

$$I^{\alpha,\psi}_{a+}f(t)=\frac{\Gamma(\delta)}{\Gamma(\alpha+\delta)}(\psi(t)-\psi(a))^{\alpha+\delta-1}.$$


\item $\displaystyle I^{\alpha,\psi}_{a+}{^{C}D^{\alpha,\psi}_{a+}}[f(x,t)]=f(x,t)-
\sum_{k=0}^{n-1}\frac{f^{[k],\psi}(x,a)}{k!}(\psi(t)-\psi(a))^k,$ where\\ $n-1<\alpha<n$ with
$n\in\mathbb{N}.$
\end{enumerate}

\end{property}
\begin{definition}
Let $\alpha>0$ and $a>0$. The one-parameter Mittag-Leffler function has the power series representation \textnormal{\cite{almeida,mittag}}
\begin{eqnarray}
E_{\alpha}[\psi(t)-\psi(a)]=\sum_{m=1}^{\infty}\frac{(\psi(t)-\psi(a))^m}{\Gamma(m\alpha+1)}. \label{ML}
\end{eqnarray}
\end{definition}

\section{Homotopy analysis method}
\label{sec:3}

In this section we introduce the basic ideas of the HAM by means of the 
description of general nonlinear problems.\newpage

We consider the following nonlinear differential equation in a general form
\begin{eqnarray}
\mathcal{N}[u(x,t)]=0, \label{eq-nl}
\end{eqnarray}
where $\mathcal{N}$ is a nonlinear differential operator, $x$ and $t$ are independent variables and $u$ is
an unknown function. We then construct the so-called zero-order deformation equation
\begin{eqnarray}
(1-p)\mathcal{L}[\varphi(x,t;p)-u_0(x,t)]=phH(x,t)\mathcal{N}[\phi(x,t;p)], \label{ordem-zero}
\end{eqnarray}
where $p\in[0,1]$ is an embedding parameter, $\hbar\neq{0}$ is an auxiliary parameter, $H(x,t)$ is an
auxiliary function and $\phi(x,t;p)$ is a function of $x$, $t$ and $p$. Let $u_0(x,t)$ be an 
initial approximation of \textnormal{Eq.(\ref{eq-nl})} and $\mathcal{L}={^{C}D_{a+}^{\alpha,\psi}}$ denotes 
an auxiliary linear differential operator with the property
$$\mathcal{L}[\phi(x,t)]=0, \quad\quad {\mbox{for}} \quad\quad \phi(x,t)=0.$$
When $p=0$ and $p=1$, we have
$$\phi(x,t;0)=u_0(x,t), \qquad {\mbox{and}} \qquad \phi(x,t;1)=u(x,t),$$
respectively. As the embedding parameter $p$ increases from $0$ to $1$, the solution $\phi(x,t;p)$ 
depends upon the embedding parameter $p$ and varies from the initial approximation $u_0(x,t)$ 
to the solution $u(x,t)$. 

Expanding $\phi(x,t;p)$ in a Taylor's series with respect to $p$, we have
\begin{eqnarray}
\phi(x,t;p)=u_0(x,t)+\sum_{m=1}^{\infty}u_m(x,t)p^{m}, \label{serie}
\end{eqnarray}
where
$$u_m(x,t)=\frac{1}{m!}\frac{\partial^m}{\partial p^m}\phi(x,t;p)\biggl|_{p=0}.$$
Assume that the auxiliary parameter $\hbar$, the auxiliary function $H(x,t)$, the initial 
approximation $u_0(x,t)$, and the auxiliary linear operator $\mathcal{L}={^{C}D_{a+}^{\alpha,\psi}}$ 
are so properly chosen that the series, \textnormal{Eq.(\ref{serie})}, converges at $p=1$. Then, 
the series \textnormal{Eq.(\ref{serie})}, at $p=1$, becomes
$$u(x,t)=\phi(x,t;1)=u_m(x,t)=u_0(x,t)+\sum_{m=1}^{\infty}u_m(x,t).$$
Differentiating \textnormal{Eq.(\ref{ordem-zero})}, $m$ times with respect to $p$, then setting $p=0$, 
and dividing it by $m!$, we obtain the $m$th-order deformation equation
\begin{eqnarray}
\mathcal{L}[u_m(x,t)-\mathcal{X}_m u_{m-1}(x,t)]=\hbar H(x,t)R_m(\vec{u}_{m-1},x,t), \label{mth}
\end{eqnarray}
with $\vec{u}_n=\{u_0(x,t),u_1(x,t),\ldots,u_n(x,t)\}$ and 
$$R_m(\vec{u}_{m-1},x,t)=\frac{1}{(m-1)!}\frac{\partial^{m-1}}{\partial p^{m-1}}\mathcal{N}[\phi(x,t;p)]\biggl|_{p=0}$$
where we have introduced the notation
\begin{eqnarray}
\mathcal{X}_{m}=\left\{
\begin{array}{lcl}
0, \quad m\leq{1},\\
1, \quad m>1. \label{x_m}
\end{array}\right.
\end{eqnarray}
Operating the fractional integral operator $I_{a+}^{\alpha,\psi}$, given by \textnormal{Eq.(\ref{integral})}, on both sides 
of \textnormal{Eq.(\ref{mth})}, we have
\begin{eqnarray}
u_m(x,t)&=&\mathcal{X}_{m}u_{m-1}(x,t)-\mathcal{X}_{m}\sum_{k=0}^{n-1}\frac{u_{m-1}^{[k],\psi}(x,a)}{k!}(\psi(t)-\psi(a))^k
\nonumber\\
&+&\hbar H(x,t)I^{\alpha,\psi}_{a+}[R_m(\vec{u}_{m-1},x,t)], \quad\quad m\geq{1}. \label{sol-mth}
\end{eqnarray}
Thus, we obtain $u_1(x,t), u_2(x,t),\cdots$ by means of \textnormal{Eq.(\ref{sol-mth})}. 
The $M$th-order approximation of $u(x,t)$ is given by
$$
u(x,t)=\sum_{m=0}^{M}u_m(x,t),
$$
and for $M\rightarrow\infty$, we get an accurate approximation of \textnormal{Eq.(\ref{eq-nl})}.
\section{Applications} 
\label{sec:4}

In this section we apply the HAM to solving linear and nonlinear fractional
partial differential equations.\\

\begin{application}
Let $t>0$, $x>0$ and $u=u(x,t)$. Consider the linear time-fractional diffusion equation 
\textnormal{\cite{jafari2006,jafari2009}}
\begin{eqnarray}
{^{C}{D}^{\alpha,\psi}_{a+}}u=\frac{\partial^2 u}{\partial x^2}+u, \quad\quad 0<\alpha<1, \label{dif-eq}
\end{eqnarray}
whose solution satisfies the initial condition
\begin{eqnarray}
u(x,a)=\cos(\pi x). \label{cond-inic}
\end{eqnarray}
In order to solve \textnormal{Eq.(\ref{dif-eq})} by means of HAM, satisfying the initial condition given by
\textnormal{Eq.(\ref{cond-inic})}, it is convenient to choose the initial approximation 
\begin{eqnarray}
u_0(x,t)=\cos(\pi x) \label{init-approx}
\end{eqnarray}
and the linear differential operator 
\begin{eqnarray*}
\mathcal{L}[\phi(x,t;p)]={^{C}}{D}^{\alpha,\psi}_{a+}[\phi(x,t;p)],
\end{eqnarray*}
satisfying the property
$$\mathcal{L}[c]=0,$$
where $c$ is an arbitrary constant. We define the nonlinear differential operator
\begin{eqnarray}
\mathcal{N}[\phi(x,t;p)]={^{C}}{D}^{\alpha,\psi}_{a+}[\phi(x,t;p)]-\frac{\partial^2}{\partial x^2}[\phi(x,t;p)]-
\phi(x,t;p). \label{op-nl}
\end{eqnarray}
Using \textnormal{Eq.(\ref{op-nl})} and the assumption $H(x,t)=1$ we construct the zero-order deformation equation
\begin{eqnarray}
(1-p)\mathcal{L}[\phi(x,t;p)-u_0(x,t)]=p\hbar\mathcal{N}[\phi(x,t;p)].
\end{eqnarray}
Obviously, when $p=0$ and $p=1$, we get
\begin{eqnarray*}
\phi(x,t;0)=u_0(x,t) \quad\quad \mbox{and} \quad\quad \phi(x,t;1)=u(x,t),
\end{eqnarray*}
respectively. So the $m$th-order deformation equation is
\begin{eqnarray}
\mathcal{L}[u_m(x,t)-\mathcal{X}_{m}u_{m-1}(x,t)]=\hbar R_{m}(\vec{u}_{m-1},x,t), \label{ordem-m}
\end{eqnarray}
subject to the initial condition $u_m(x,a)=0$ where $\mathcal{X}_m$ is defined by \textnormal{Eq.(\ref{x_m})} and
$$R_{m}(\vec{u}_{m-1},x,t)={^{C}}{D}^{\alpha,\psi}_{a+}u_{m-1}(x,t)-\frac{\partial^2}{\partial x^2}u_{m-1}(x,t)-
u_{m-1}(x,t).$$
Now we apply the integral fractional operator $I_{a+}^{\alpha,\psi}$ on both sides of \textnormal{Eq.(\ref{ordem-m})} to get
\begin{eqnarray*}
{^{C}}{D}^{\alpha,\psi}_{a+}[u_m(x,t)-\mathcal{X}_{m}u_{m-1}(x,t)]=\hbar I^{\alpha,\psi}_{a+}
\left[{^{C}}{D}^{\alpha,\psi}_{a+}u_{m-1}(x,t)-\frac{\partial^2}{\partial x^2}u_{m-1}(x,t)-
u_{m-1}(x,t)\right],
\end{eqnarray*}
whose solution has the form
\begin{eqnarray*}
u_m(x,t)&-&\sum_{k=0}^{n-1}\frac{u_{m}^{[k],\psi}(x,a)}{k!}
(\psi(t)-\psi(a))^k-\mathcal{X}_{m}u_{m-1}(x,t)\\
&+&\mathcal{X}_{m}\sum_{k=0}^{n-1}\frac{u_{m-1}^{[k],\psi}(x,a)}{k!}
(\psi(t)-\psi(a))^k=\hbar\left[u_{m-1}(x,t)-\sum_{k=0}^{n-1}\frac{u_{m-1}^{[k],\psi}(x,a)}{k!}
(\psi(t)-\psi(a))^k\right.\\
&-&\left. I^{\alpha,\psi}_{a+}\left(\frac{\partial^2}{\partial x^2}u_{m-1}(x,t)+u_{m-1}(x,t)\right)\right],
\quad\quad m\geq{1}.
\end{eqnarray*}
For $0<\alpha<1$, then $n=1$, we can rewrite the above equation as
\begin{eqnarray}
u_m(x,t)&=&(\mathcal{X}_m+\hbar)u_{m-1}(x,t)-(\mathcal{X}_{m}+\hbar)u_{m-1}(x,a)\nonumber\\
&-&\hbar I^{\alpha,\psi}_{a+}
\left[\frac{\partial^2}{\partial x^2}u_{m-1}(x,t)+u_{m-1}(x,t)\right], \quad\quad m\geq{1}. \label{u_m}
\end{eqnarray}
From \textnormal{Eq.(\ref{init-approx})} and \textnormal{Eq.(\ref{u_m})}, we obtain
\begin{eqnarray*}
u_0(x,t)&=&\cos(\pi x),\\
u_1(x,t)&=&-\hbar I^{\alpha,\psi}_{a+}\left[\frac{\partial^2}{\partial x^2}u_0(x,t)+u_0(x,t)\right]=
-\hbar(1-\pi^2)\cos(\pi x)\frac{(\psi(t)-\psi(a))^{\alpha}}{\Gamma(\alpha+1)},\\
u_2(x,t)&=&-(1+\hbar)\hbar(1-\pi^2)\cos(\pi x)\frac{(\psi(t)-\psi(a))^{\alpha}}{\Gamma(\alpha+1)}+
\hbar^2(1-\pi^2)^2\cos(\pi x)\frac{(\psi(t)-\psi(a))^{2\alpha}}{\Gamma(2\alpha+1)},\\
&\vdots&
\end{eqnarray*}
An accurate approximation of \textnormal{Eq.(\ref{dif-eq})} is given by
$$u(x,t)=u_0(x,t)+u_1(x,t)+u_2(x,t)+\cdots,$$
and, when $\hbar=-1$, we have
\begin{eqnarray*}
u(x,t)&=&\cos(\pi x)\left[1+\frac{(1-\pi^2)}{\Gamma(\alpha+1)}(\psi(t)-\psi(a))^{\alpha}+
\frac{[(1-\pi^2)](\psi(t)-\psi(a))^{\alpha}]^{2}}{\Gamma(2\alpha+1)}\right.\\
&+&\left.\frac{[(1-\pi^2)(\psi(t)-\psi(a))^{\alpha}]^{3}}{\Gamma(3\alpha+1)}+\cdots\right]\\
&=&\cos(\pi x)\left[1+\sum_{m=1}^{\infty}\frac{[(1-\pi^2)(\psi(t)-\psi(a))^{\alpha}]^{m}}{\Gamma(m\alpha+1)}\right],
\end{eqnarray*}
or in terms of the one-parameter Mittag-Lefler function \textnormal{Eq.(\ref{ML})},
\begin{eqnarray}
u(x,t)=\cos(\pi x)E_{\alpha}[(1-\pi^2)(\psi(t)-\psi(a))^{\alpha}]. \label{sol}
\end{eqnarray}
We note two important special cases of \textnormal{Eq.(\ref{sol})}. First taking $\psi(t)=t$ and $a=0$. In this case
the solution \textnormal{Eq.(\ref{sol})} takes the form
\begin{eqnarray}
u(x,t)=\cos(\pi x)\,E_{\alpha}[(1-\pi^2)t^{\alpha}] . \label{sol-t} 
\end{eqnarray}
\textnormal{Eq.(\ref{sol-t})} recovers the solutions found by Jafari and Seifi \textnormal{\cite{jafari2009}} obtained by means of the
HAM and Jafari and Daftardar-Gejji \textnormal{\cite{jafari2006}} using the Adomian
decomposition method. 

On the other hand, if $\psi(t)=\ln t$ and $a>0$, the solution \textnormal{Eq.(\ref{sol})} becomes
\begin{eqnarray}
u(x,t)=\cos(\pi x)\,E_{\alpha}\left[(1-\pi^2)\left(\ln\frac{t}{a}\right)^{\alpha}\right]. \label{sol-lnt}
\end{eqnarray}
\begin{figure}[H]
\centering
\includegraphics[width=0.65\textwidth]{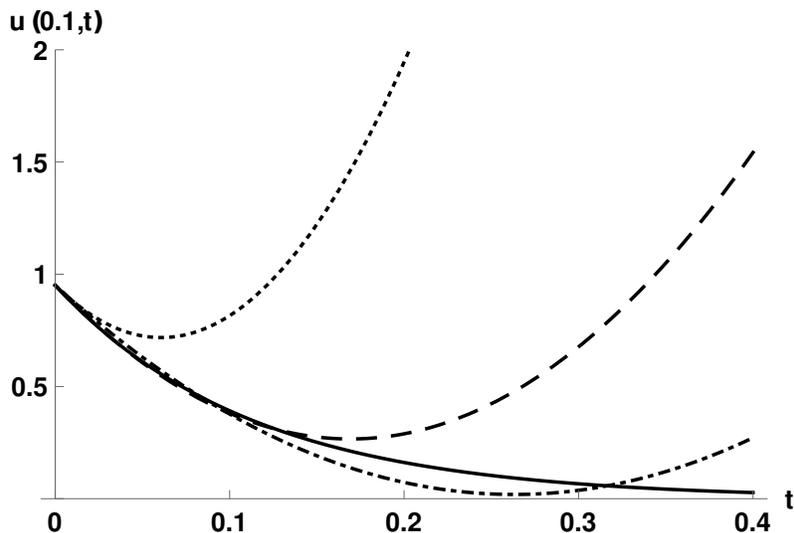}
\caption{Approximate solutions $u(0.1,t)$ using 3-terms and exact solution of \textnormal{Eq.(\ref{dif-eq})} subject the initial condition
\textnormal{Eq.(\ref{cond-inic})} with $\psi(t)=t$, $a=0$ and $\alpha\rightarrow{1}$. 
Solid line (exact solution): $\hbar=-1$, dashdotted: $\hbar=-0.6$, dashed: $\hbar=-0.8$, and dotted: $\hbar=-1.3$.}
\end{figure}
\begin{figure}[H]
\centering
\includegraphics[width=0.65\textwidth]{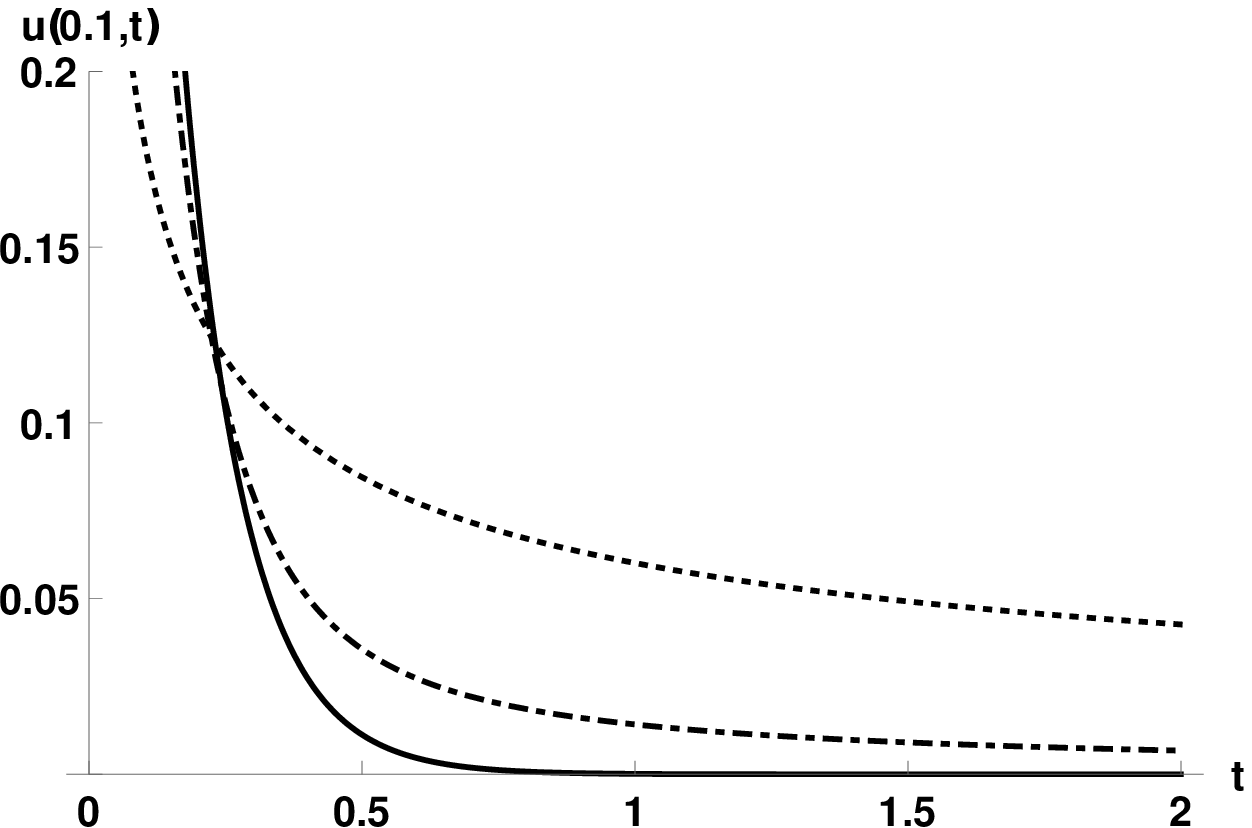}
\caption{Exact solutions $u(0.1,t)$ using \textnormal{Eq.(\ref{sol-t})}. Solid line: $\alpha\rightarrow{1}$,
dashdotted: $\alpha=0.9$, and dotted: $\alpha=0.5$.}
\end{figure}
\begin{figure}[H]
\centering
\includegraphics[width=0.65\textwidth]{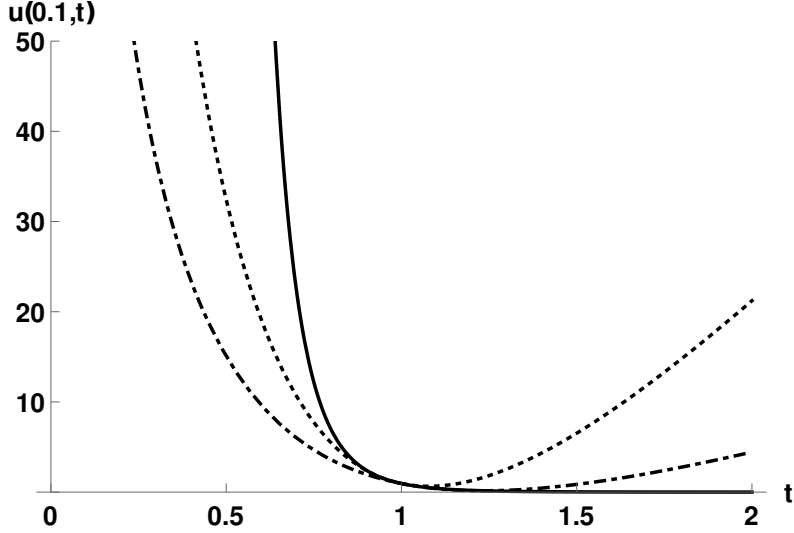}
\caption{Approximate solutions $u(0.1,t)$ using 3-terms, $\psi(t)=\ln t$, $a=1$ and $\alpha\rightarrow{1}$. Solid line 
(exact solution using \textnormal{Eq.(\ref{sol-lnt}))}: $\hbar=-1$, dashdotted: $\hbar=-0.7$, and dotted: $\hbar=-1.2$.}
\end{figure}
\end{application}
\begin{application}

Let $t>0$, $x >0$ and $u=u(x,t)$. Consider the nonlinear time-fractional gas-dynamic equation \textnormal{\cite{shone}}
\begin{eqnarray}
{^{C}{D}^{\alpha,\psi}_{a+}}u+u\cdot\frac{\partial u}{\partial x}-u+u^2=0, \quad\quad 0<\alpha<1 \label{dif-eq2}
\end{eqnarray}
whose solution satisfies the initial condition
\begin{eqnarray}
u(x,a)=\mbox{e}^{-x}. \label{cond-inic-2}
\end{eqnarray}
In order to solve \textnormal{Eq.(\ref{dif-eq2})}, we choose the initial approximation
\begin{eqnarray*}
u_0(x,t)=\mbox{e}^{-x} \label{init-approx-2}
\end{eqnarray*}
and the linear operator
\begin{eqnarray*}
\mathcal{L}[\phi(x,t;p)]={^{C}}{D}^{\alpha,\psi}_{a+}[\phi(x,t;p)],
\end{eqnarray*}
with the property $\mathcal{L}[c]=0$, where $c$ is a  constant.
From \textnormal{Eq.(\ref{dif-eq2})}, we define the nonlinear differential operator
\begin{eqnarray*}
\mathcal{N}[\phi(x,t;p)]={^{C}}{D}^{\alpha,\psi}_{a+}[\phi(x,t;p)]+\phi(x,t;p)\cdot\frac{\partial}{\partial x}[\phi(x,t;p)]-
\phi(x,t;p)+[\phi(x,t;p)]^2.
\end{eqnarray*}
Taking $H(x,t)=1$, we construct the zero-order deformation equation
\begin{eqnarray}
(1-p)\mathcal{L}[\phi(x,t;p)-u_0(x,t)]=p\hbar\mathcal{N}[\phi(x,t;p)].
\end{eqnarray}
Obviously, when $p=0$ and $p=1$, we get
\begin{eqnarray*}
\phi(x,t;0)=u_0(x,t)=\mbox{e}^{-x} \quad\quad \mbox{and} \quad\quad \phi(x,t;1)=u(x,t),
\end{eqnarray*}
respectively. The $m$th-order deformation equation is given by
\begin{eqnarray}
\mathcal{L}[u_m(x,t)-\mathcal{X}_{m}u_{m-1}(x,t)]=\hbar R_{m}(\vec{u}_{m-1},x,t), \label{ordem-m2}
\end{eqnarray}
subject to the initial condition $u_m(x,a)=0$, where
\begin{eqnarray}
R_{m}(\vec{u}_{m-1},x,t)&=&{^{C}}{D}^{\alpha,\psi}_{a+}u_{m-1}(x,t)+\sum_{i=0}^{m-1}u_i(x,t)\cdot\frac{\partial}{\partial x}u_{m-1-i}(x,t)-u_{m-1}(x,t)\nonumber\\
&+&\sum_{i=0}^{m-1}u_i(x,t)u_{m-1-i}(x,t).
\end{eqnarray}
Operating the fractional integral operator $I_{a+}^{\alpha,\psi}$ on both sides of \textnormal{Eq.(\ref{ordem-m2})},
we have
\begin{eqnarray}
&&u_m(x,t)=(\mathcal{X}_m+\hbar)u_{m-1}(x,t)-(\mathcal{X}_{m}+\hbar)u_{m-1}(x,a) \label{u_m-2}\\
&&+\hbar I^{\alpha,\psi}_{a+}
\left[\sum_{i=0}^{m-1}u_i(x,t)\cdot\frac{\partial}{\partial x}u_{m-1-i}(x,t)-u_{m-1}(x,t)+
\sum_{i=0}^{m-1}u_i(x,t)u_{m-1-i}(x,t)\right]. \nonumber
\end{eqnarray}
In this way, we obtain
\begin{eqnarray*}
u_0(x,t)&=&\mbox{e}^{-x},\\
u_1(x,t)&=&\hbar I^{\alpha,\psi}_{a+}\left\{u_0(x,t)\cdot\frac{\partial}{\partial x}u_0(x,t)-u_0(x,t)+[u_0(x,t)]^2]\right\}
=-\hbar\mbox{e}^{-x}\frac{(\psi(t)-\psi(a))^{\alpha}}{\Gamma(\alpha+1)},\\
u_2(x,t)&=&(1+\hbar)u_1(x,t)-(1+\hbar)u_1(x,a)+\hbar{I^{\alpha,\psi}_{a+}}\left\{u_0(x,t)\cdot\frac{\partial}{\partial x}u_1(x,t)\right.\\
&+&\left. u_1(x,t)\cdot\frac{\partial}{\partial x}u_0(x,t)-u_1(x,t)+2u_1(x,t)\cdot u_0(x,t)\right\}\\
&=&-(1+\hbar)\hbar\,\mbox{e}^{-x}\frac{(\psi(t)-\psi(a))^{\alpha}}{\Gamma(\alpha+1)}+
{\hbar^2}\mbox{e}^{-x}\frac{(\psi(t)-\psi(a))^{2\alpha}}{\Gamma(2\alpha+1)},\\
u_3(x,t)&=&(1+\hbar)u_2(x,t)-(1+\hbar)u_2(x,a)+\hbar{I^{\alpha,\psi}_{a+}}\left\{u_0(x,t)\cdot\frac{\partial}{\partial x}u_2(x,t)\right.\\
&+&\left. u_1(x,t)\cdot\frac{\partial}{\partial x}u_1(x,t)+u_2(x,t)\cdot\frac{\partial}{\partial x}u_0(x,t)
-u_2(x,t)\right.\\
&+&2u_0(x,t)\cdot u_2(x,t)+[u_1(x,t)]^2\bigg\}\\
&=&-(1+\hbar)^2 \hbar\,\mbox{e}^{-x}\frac{(\psi(t)-\psi(a))^{\alpha}}{\Gamma(\alpha+1)}+
2{\hbar^2}(1+\hbar)\mbox{e}^{-x}\frac{(\psi(t)-\psi(a))^{2\alpha}}{\Gamma(2\alpha+1)}\\
&-&{\hbar^3}\mbox{e}^{-x}\frac{(\psi(t)-\psi(a))^{3\alpha}}{\Gamma(3\alpha+1)},\\
&\vdots&
\end{eqnarray*}
The solution $u(x,t)$ is given by
$$u(x,t)=u_0(x,t)+u_1(x,t)+u_2(x,t)+\cdots$$
and if $\hbar=-1$, we have
\begin{eqnarray}
u(x,t)&=&\mbox{e}^{-x}\left[1+\frac{(\psi(t)-\psi(a))^{\alpha}}{\Gamma(\alpha+1)}+
\frac{(\psi(t)-\psi(a))^{2\alpha}}{\Gamma(2\alpha+1)}+
\frac{(\psi(t)-\psi(a))^{3\alpha}}{\Gamma(3\alpha+1)}+\cdots\right]\nonumber\\
&=&\mbox{e}^{-x}E_{\alpha}[(\psi(t)-\psi(a))^{\alpha}]. \label{sol-2}
\end{eqnarray}
In particular, if $\psi(t)=t$ and $a=0$, we can write the solution as
\begin{eqnarray}
u(x,t)=\mbox{e}^{-x}E_{\alpha}(t^{\alpha}) \label{sol-t-2}
\end{eqnarray}
and if $\alpha\rightarrow{1}$, we obtain the solution found by Shone and Patra \textnormal{\cite{shone}}
using the fractional complex transform and a new iterative method, this is, $u(x,t)=e^{t-x}$. On the other hand, if
$\psi(t)=\ln t$ and $a>0$, we obtain
\begin{eqnarray}
u(x,t)=\mbox{e}^{-x}\,E_{\alpha}\left[\left(\ln\frac{t}{a}\right)^{\alpha}\right]. \label{sol-lnt-2}
\end{eqnarray}
\begin{figure}[H]
\centering
\includegraphics[width=0.65\textwidth]{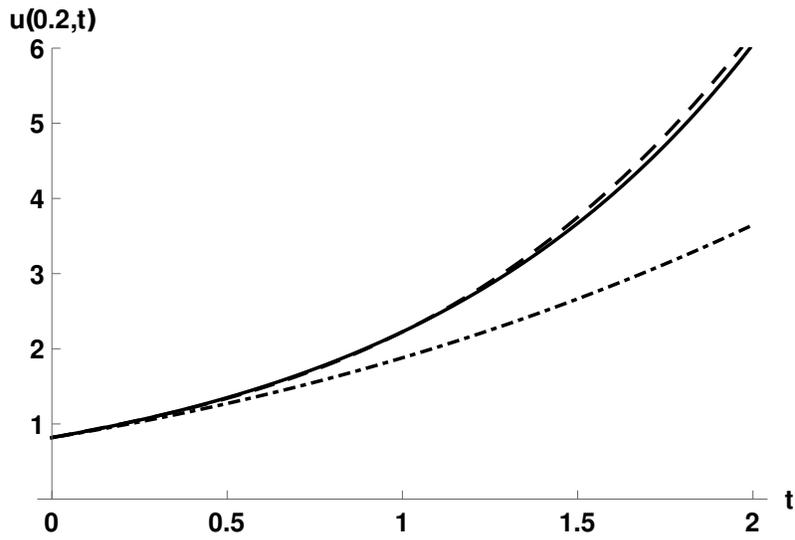}
\caption{Approximate solutions $u(0.2,t)$ using 4-terms and exact solution of \textnormal{Eq.(\ref{dif-eq2})} subject the initial
condition \textnormal{Eq.(\ref{cond-inic-2})} with $\psi(t)=t$, $a=0$ and $\alpha\rightarrow{1}$. Solid line (exact solution): $\hbar=-1$,
dashdotted: $\hbar=-0.6$, and dashed: $\hbar=-1.4$.}
\end{figure}
\begin{figure}[H]
\centering
\includegraphics[width=0.65\textwidth]{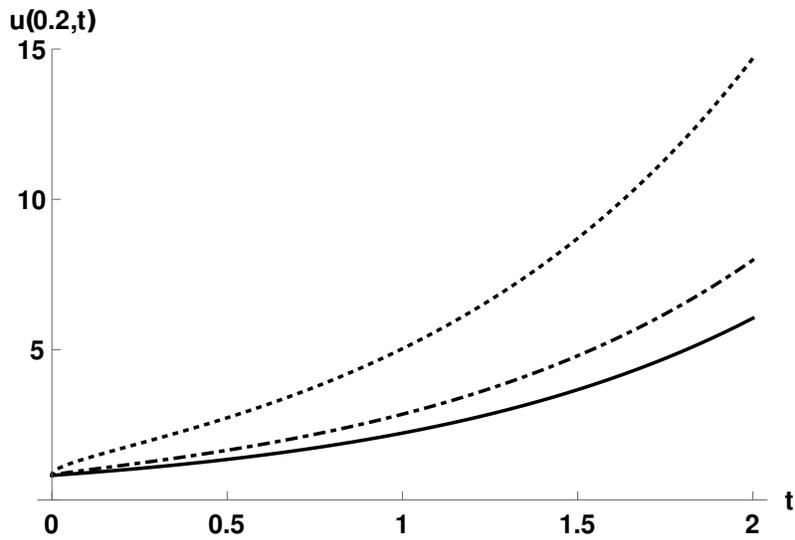}
\caption{Exact solutions $u(0.2,t)$ using \textnormal{Eq.(\ref{sol-t-2})}. Solid line: $\alpha\rightarrow{1}$,
dashdotted: $\alpha=0.75$, and dotted: $\alpha=0.4$.}
\end{figure}
\begin{figure}[H]
\centering
\includegraphics[width=0.65\textwidth]{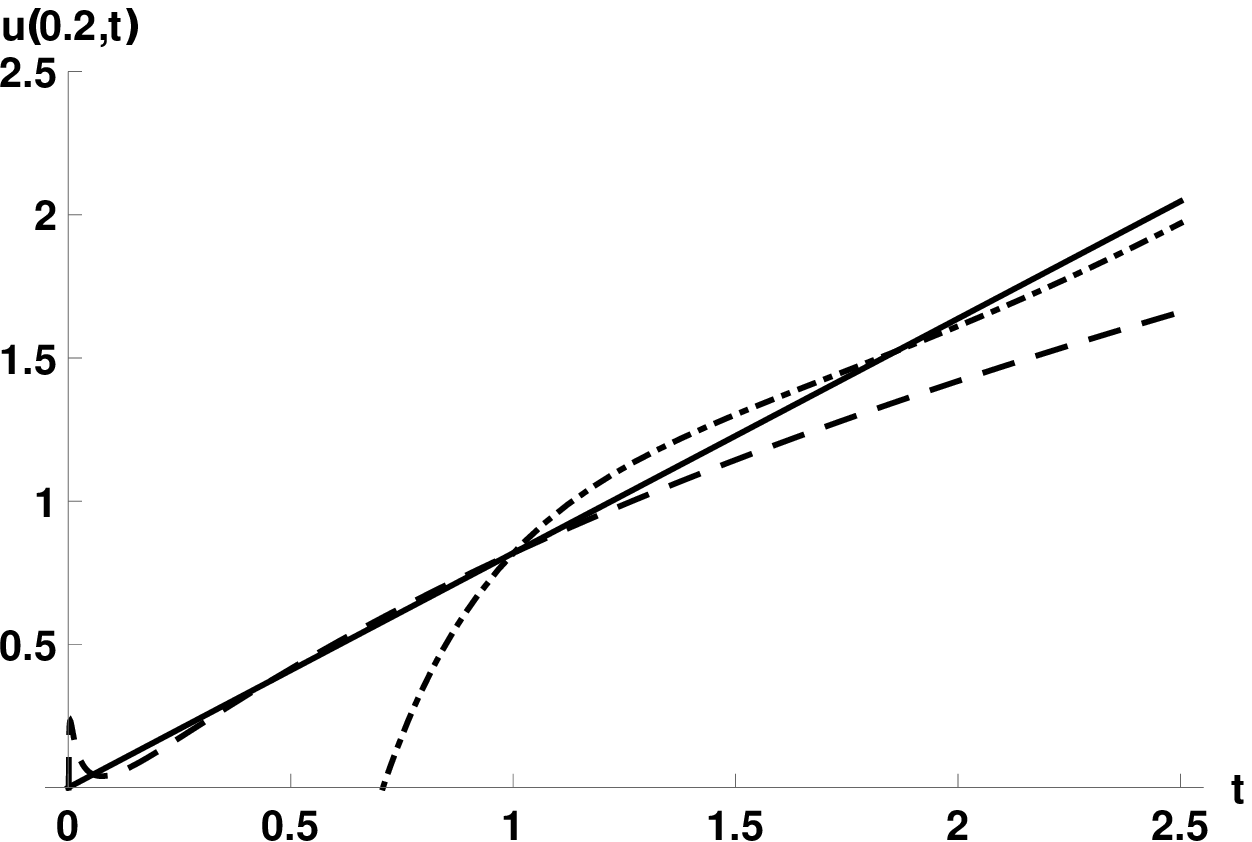}
\caption{Approximate solutions $u(0.2,t)$ using 4-terms, $\psi(t)=\ln t$, $a=1$ and $\alpha\rightarrow{1}$. 
Solid line (exact solution using \textnormal{Eq.(\ref{sol-lnt-2}))}: $\hbar=-1$, dashdotted: $\hbar=-2$, and dashed: $\hbar=-0.5$.}
\end{figure}
\end{application}
\begin{application}
\noindent Let $t>0$, $x>0$, $u=u(x,t)$ and $0<\alpha<1$. Consider the following nonlinear time-fractional KdV equation 
\textnormal{\cite{momani,rehman}}
\begin{eqnarray}
{^{C}{D}^{\alpha,\psi}_{a+}}u(x,t)-\frac{\partial}{\partial x}[u(x,t)]^2+\frac{\partial}{\partial x}
\left[u(x,t)\frac{\partial^2}{\partial x^2}u(x,t)\right]=0,  \label{dif-eq3}
\end{eqnarray}
whose solution satisfies the initial condition
\begin{eqnarray}
u(x,a)=\sinh^2 \left(\frac{x}{2}\right). \label{cond-inic-3}
\end{eqnarray}
Let $u_0(x,t)$ denote an initial approximation of $u(x,t)$, this is,
\begin{eqnarray}
u_0(x,t)=\sinh^2 \left(\frac{x}{2}\right) \label{init-approx-3}
\end{eqnarray}
and we choose the linear differential operator $\mathcal{L}={^{C}{D}^{\alpha,\psi}_{a+}}$, with the condition
$\mathcal{L}[c]=0$ where $c$ is a constant.
From \textnormal{Eq.(\ref{dif-eq3})}, we define the nonlinear differential operator
\begin{eqnarray*}
\mathcal{N}[\phi(x,t;p)]={^{C}}{D}^{\alpha,\psi}_{a+}[\phi(x,t;p)]+\phi(x,t;p)\cdot\frac{\partial}{\partial x}[\phi(x,t;p)]-
\phi(x,t;p)+[\phi(x,t;p)]^2.
\end{eqnarray*}
With the choice $H(x,t)=1$ we have the zero-order deformation equation
\begin{eqnarray}
(1-p)\mathcal{L}[\phi(x,t;p)-u_0(x,t)]=p\hbar\mathcal{N}[\phi(x,t;p)].
\end{eqnarray}
Obviously, when $p=0$ and $p=1$, we get
\begin{eqnarray*}
\phi(x,t;0)=u_0(x,t)=\sinh^2 \left(\frac{x}{2}\right) \quad\quad \mbox{and} \quad\quad \phi(x,t;1)=u(x,t),
\end{eqnarray*}
respectively. The $m$th-order deformation equation can be expressed by
\begin{eqnarray}
\mathcal{L}[u_m(x,t)-\mathcal{X}_{m}u_{m-1}(x,t)]=\hbar R_{m}(\vec{u}_{m-1},x,t), \label{ordem-m3}
\end{eqnarray}
where
\begin{eqnarray*}
R_{m}(\vec{u}_{m-1},x,t)&=&{^{C}}{D}^{\alpha,\psi}_{a+}u_{m-1}(x,t)-\frac{\partial}{\partial x}\left[\sum_{i=0}^{m-1}u_i(x,t)\cdot u_{m-1-i}(x,t)\right.\\
&-&\left.\sum_{i=0}^{m-1}u_i(x,t)\frac{\partial^2}{\partial x^2}u_{m-1-i}(x,t)\right].
\end{eqnarray*}
Applying the fractional operator $I_{a+}^{\alpha,\psi}$ to this equation we find
\begin{eqnarray}
u_m(x,t)&=&(\mathcal{X}_m+\hbar)u_{m-1}(x,t)-(\mathcal{X}_{m}+\hbar)u_{m-1}(x,a) \label{u_m-3}\\
&-&\hbar I^{\alpha,\psi}_{a+}
\left[\frac{\partial}{\partial x}\left(\sum_{i=0}^{m-1}u_i(x,t)\cdot u_{m-1-i}(x,t)-\sum_{i=0}^{m-1}u_{i}(x,t)
\frac{\partial^2}{\partial x^2}u_{m-1-i}(x,t)\right)\right]. \nonumber
\end{eqnarray}
Thereafter, we successively obtain
\begin{eqnarray*}
u_0(x,t)&=&\sinh^2 \left(\frac{x}{2}\right),\\
u_1(x,t)&=&-\hbar I^{\alpha,\psi}_{a+}\left[\frac{\partial}{\partial x}\left([u_0(x,t)]^2-
u_0(x,t)\frac{\partial^2}{\partial x^2}u_0(x,t)\right)\right]=\hbar\sinh(x)\frac{(\psi(t)-\psi(a))^{\alpha}}{4\Gamma(\alpha+1)},\\
u_2(x,t)&=&(1+\hbar)u_1(x,t)-\hbar{I^{\alpha,\psi}_{a+}}\left[\frac{\partial}{\partial x}\left(2u_0(x,t)\cdot u_1(x,t)-
u_0(x,t)\frac{\partial^2}{\partial x^2}u_1(x,t)-u_1(x,t)\frac{\partial^2}{\partial x^2}u_0(x,t)\right)\right]\\
&=&(1+\hbar)\hbar\,\sinh(x)\,\frac{(\psi(t)-\psi(a))^{\alpha}}{4\Gamma(\alpha+1)}+
\hbar^2\cosh(x)\,\frac{(\psi(t)-\psi(a))^{2\alpha}}{8\Gamma(2\alpha+1)},\\
&\vdots&
\end{eqnarray*}
The second-order approximation of $u(x,t)$ is
\begin{eqnarray*}
u(x,t)&=&\sinh^2 \left(\frac{x}{2}\right)+\hbar\sinh(x)\frac{(\psi(t)-\psi(a))^{\alpha}}{4\Gamma(\alpha+1)}+
(1+\hbar)\hbar\sinh(x)\,\frac{(\psi(t)-\psi(a))^{\alpha}}{4\Gamma(\alpha+1)}\\
&+&\hbar^2\cosh(x)\,\frac{(\psi(t)-\psi(a))^{2\alpha}}{8\Gamma(2\alpha+1)}.
\end{eqnarray*}
Taking $\hbar=-1$, we have
\begin{eqnarray}
u(x,t)=\sinh^2 \left(\frac{x}{2}\right)-\frac{\sinh(x)}{4}\frac{(\psi(t)-\psi(a))^{\alpha}}{\Gamma(\alpha+1)}+
\frac{\cosh(x)}{8}\frac{(\psi(t)-\psi(a))^{2\alpha}}{\Gamma(2\alpha+1)}. \label{sol-2nd}
\end{eqnarray}
If $\psi(t)=t$ and $a=0$, the second-order approximation of $u(x,t)$ \textnormal{Eq.(\ref{sol-2nd})} becomes
\begin{eqnarray}
u(x,t)=\sinh^2 \left(\frac{x}{2}\right)-\frac{t^{\alpha}}{4\Gamma(\alpha+1)}\,\sinh(x)+
\frac{t^{2\alpha}}{8\Gamma(2\alpha+1)}\,\cosh(x).\label{sol-t-3}
\end{eqnarray}
This solution is identical to the solution obtained using Rehman et al. \textnormal{\cite{rehman}}
the combination of the double Sumudu transform and homotopy perturbation method and also
obtained by Momani et al. \textnormal{\cite{momani}} by homotopy perturbation method. 

On the other hand, if $\psi(t)=\ln t$ and $a>0$, we obtain
\begin{eqnarray}
u(x,t)=\sinh^2 \left(\frac{x}{2}\right)-\frac{\sinh(x)}{4\Gamma(\alpha+1)}\left(\ln\frac{t}{a}\right)^\alpha+
\frac{\cosh(x)}{8\Gamma(2\alpha+1)}\left(\ln\frac{t}{a}\right)^{2\alpha}. \label{sol-t-4}
\end{eqnarray}
\begin{figure}[H]
\centering
\includegraphics[width=0.65\textwidth]{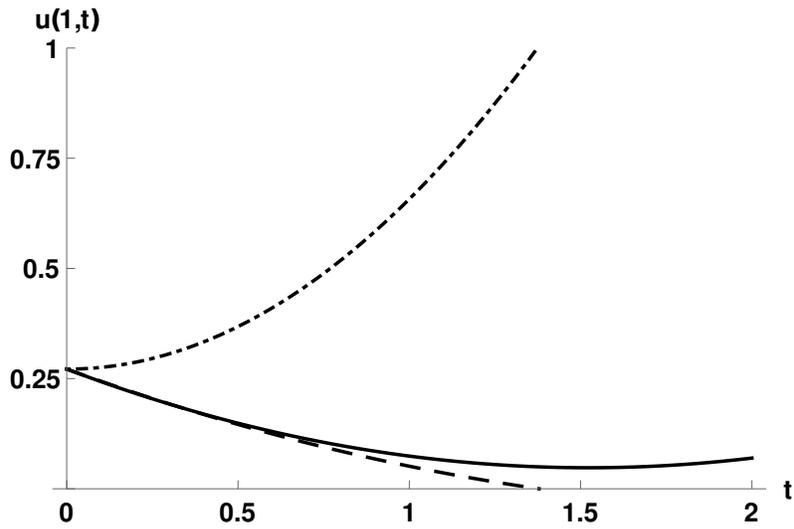}
\caption{Approximate solutions $u(1,t)$ using 3-terms, $\psi(t)=t$, $a=0$ and $\alpha\rightarrow{1}$. Solid line: $\hbar=-1$,
dashdotted: $\hbar=-2$, and dashed: $\hbar=-0.8$.}
\end{figure}
\begin{figure}[H]
\centering
\includegraphics[width=0.65\textwidth]{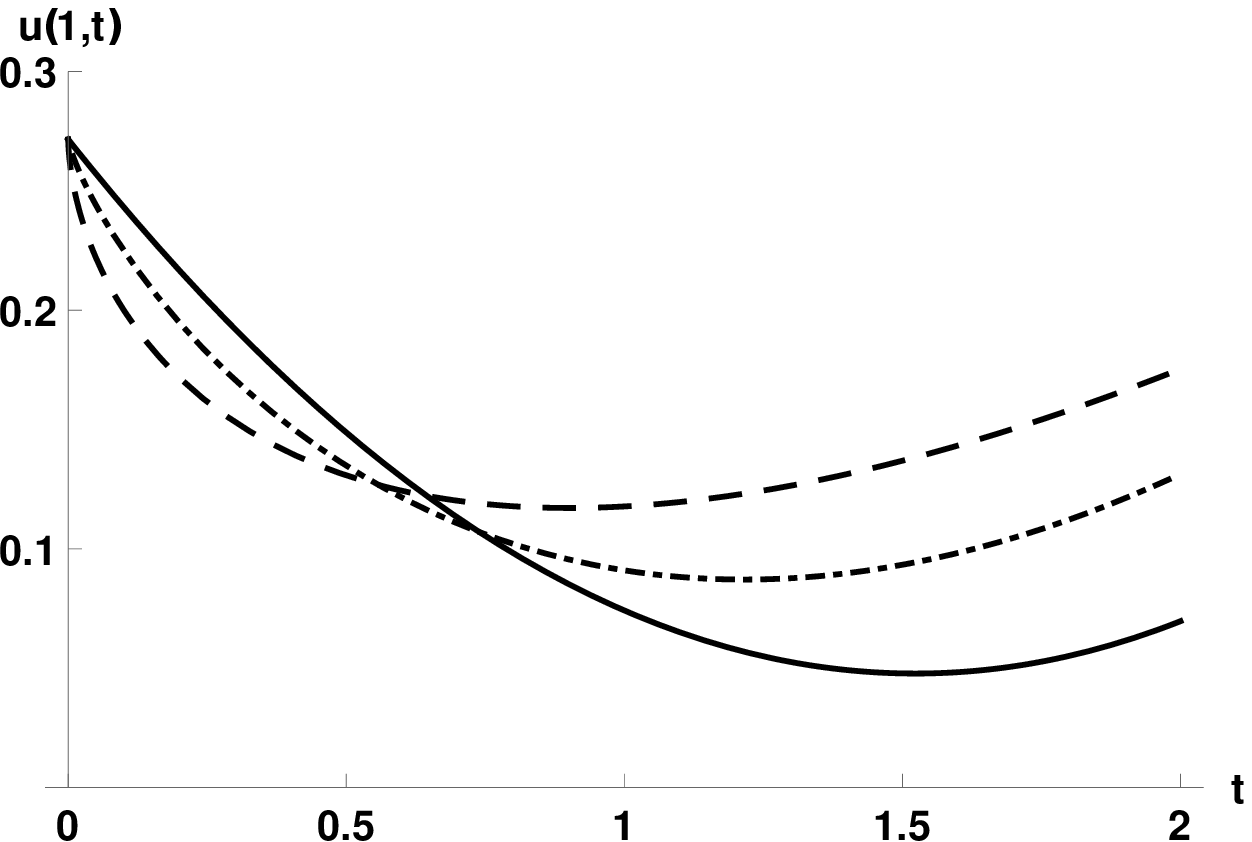}
\caption{Approximate solutions $u(1,t)$ using \textnormal{Eq.(\ref{sol-t-3})}. Solid line: $\alpha\rightarrow{1}$,
dashdotted: $\alpha=0.8$, and dashed: $\alpha=0.6$.}
\end{figure}
\begin{figure}[H]
\centering
\includegraphics[width=0.65\textwidth]{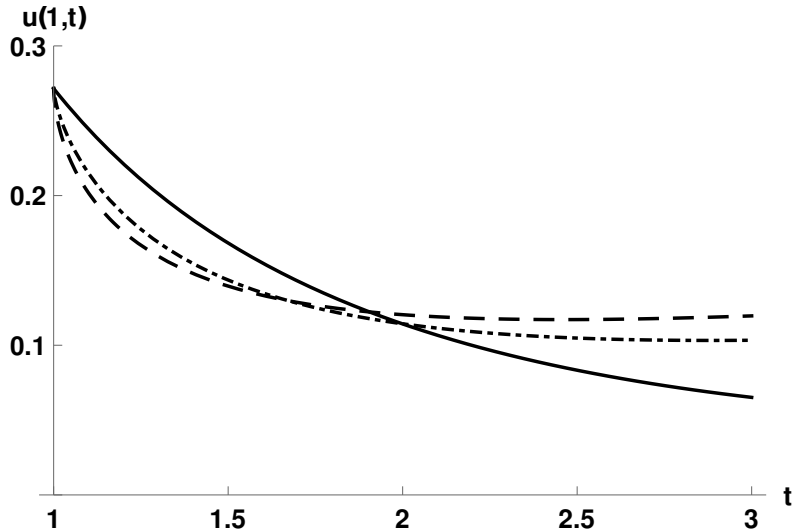}
\caption{Approximate solutions $u(1,t)$ using \textnormal{Eq.(\ref{sol-t-4})}. Solid line: $\alpha\rightarrow{1}$,
dashdotted: $\alpha=0.7$, and dashed: $\alpha=0.6$.}
\end{figure}
\end{application}
\section{Concluding remarks}
\label{sec:5}

In this paper we have presented the HAM to obtain approximate solutions for linear and nonlinear fractional
partial differential equations replacing the first order time derivative by the $\psi$-Caputo 
fractional derivative. We solve fractional partial differential equations, and to obtain
explicit series solutions, we have presented numerical solutions. Therefore,
we considered different values for $\alpha$ and for auxiliary parameter $\hbar$. It is possible to 
control the convergence region of the solution series, obtained by means of HAM, adjusting the auxiliary 
parameter $\hbar$. Mathematica has been used for draw graphs.





\begin{thebibliography}{}
%
%

\bibitem{Adomian}
G. Adomian, Solving Frontier Problems of Physics: The Decomposition Method. Kluwer Acad. Publ., Boston, (1994).

\bibitem{almeida}
R. Almeida, A Caputo fractional derivative of a function with respect to another
function, Commun. Nonlinear Sci. Numer. Simulat. 44:460--481, (2017).

\bibitem{ganjiani}
M. Ganjiani, Solution of nonlinear fractional differential equations using homotopy analysis
method, Appl. Math. Model. 34:1634--1641, (2010).

\bibitem{jafari2006}
H. Jafari and V. Daftardar-Gejji, Solving linear and nonlinear fractional diffusion and
wave equations by Adomian decomposition, Appl. Math. Comput., 180:488--497, (2006).

\bibitem{jafari2009}
H. Jafari and S. Seifi, Homotopy analysis method for solving linear and nonlinear fractional
diffusion-wave equation, Commun. Nonlinear Sci. Numer. Simulat., 14:2006--2012, (2009).

\bibitem{kilbas}
A. A. Kilbas, H. M. Srivastava and J. Trujillo, Theory and Applications of the Fractional
Differential Equations. Vol. 204. Elsevier, Amsterdam (2006).

\bibitem{Liao2}
S. J. Liao, The proposed homotopy analysis technique for the solution of nonlinear problems, Ph.D. Thesis, 
Shanghai Jiao Tong University, (1992).

\bibitem{Liao}
S. J. Liao, Beyond Perturbation: Introduction to the Homotopy Analysis Method. Chapman \& Hall, Boca Raton (2003).

\bibitem{mittag}
M. G. Mittag-Leffler, Sur la nouvelle fonction $E_{\alpha}(x)$, C. R. Acad. Sci., 137:554--558, (1903).

\bibitem{momani}
S. Momani, Z. Odibat and I. Hashim, Algorithms for nonlinear fractional partial
differential equations: a selection of numerical methods, Topol. Methods Nonlinear Anal., 
31:211--226, (2008).

\bibitem{danieco}
D. S. Oliveira and E. Capelas de Oliveira, Hilfer-Katugampola fractional derivatives, Com. Appl. Math., 
37:3672--3690, (2017).

\bibitem{rehman}
H. U. Rehman, M. S. Sallem and A. Ahmad, Combination of Homotopy Perturbation Method (HPM) and 
double Sumudu transform to solve fractional KdV equations, Open J. Math. Sci., 2:29--38, (2018).

\bibitem{graziane}
G. Sales Teodoro, J. A. Tenreiro Machado and E. Capelas de Oliveira, A review of definitions
of fractional derivatives and other operators. J. Comput. Phys., 388:195--208, (2019).

\bibitem{shone}
T. T. Shone and A. Patra, Solution for non-linear fractional partial differential equations 
using fractional complex transform, Int. J. Appl. Comput. Math., 5:90 (8 pages), (2019).

\bibitem{slota}
D. S{\l}ota, E. Hetmaniok, R. Witu{\l}a, K. Gromysz and T. Trawi\'{n}ski, 
Homotopy approach for integrodifferential equations, Math., 7:904 (15 pages), (2019).

\bibitem{vanterler}
J. Vanterler da C. Sousa and E. Capelas de Oliveira, On the $\psi$-Hilfer fractional derivative,
Commun. Nonlinear Sci. Numer. Simulat., 60:72--91, (2018).

\bibitem{Xu}
H. Xu, S. J. Liao, XC. You, Analysis of nonlinear fractional partial differential equations with the
homotopy analysis method, Commun. Nonlinear Sci. Numer. Simulat., 14:1152--1156, (2009).

\bibitem{zhang}
X. Zhang, B. Tang and Y. He, Homotopy analysis method for higher-order fractional integro-differential 
equations, Comput. Math. Appl., 62:3194--3203, (2011).

\bibitem{zurigat2010}
M. Zurigat, S. Momani, Z. Odibat and A. Alawneh, The homotopy analysis method for handling systems of fractional
differential equations, Appl. Math. Model., 34:24--35, (2010).

\bibitem{zurigat}
M. Zurigat, S. Momani and A. Alawneh, Analytical approximate solutions of systems of fractional
algebraic-differential equations by homotopy analysis method, Comput. Math. Appl.,59:1227--1235,
(2010).


\end{thebibliography}

\end{document}